\documentclass[review,hidelinks,onefignum,onetabnum]{siamart250211}

\usepackage{lipsum}
\usepackage{amsfonts}
\usepackage{graphicx}
\usepackage{epstopdf}
\usepackage{subcaption}
\usepackage{algorithmic}
\ifpdf
  \DeclareGraphicsExtensions{.eps,.pdf,.png,.jpg}
\else
  \DeclareGraphicsExtensions{.eps}
\fi

\newsiamremark{remark}{Remark}
\newsiamremark{hypothesis}{Hypothesis}
\crefname{hypothesis}{Hypothesis}{Hypotheses}
\newsiamthm{claim}{Claim}
\newsiamremark{fact}{Fact}
\crefname{fact}{Fact}{Facts}

\headers{Kernel Variational Inference Flow Filter}{Weiye Gan, Zhijun Zeng, Junqing Chen, and Zuoqiang Shi}

\title{Kernel Variational Inference Flow for Nonlinear Filtering Problem\thanks{Submitted to the editors DATE.
\funding{This work was supported by National Natural Science Foundation of China under grant 12071244.}}}

\author{Weiye Gan\thanks{Department of Mathematical Sciences, Tsinghua University, Beijing, China 
  (\email{gwy22@mails.tsinghua.edu.cn}).}\and
Zhijun Zeng\thanks{Department of Mathematical Sciences, Tsinghua University, Beijing, China 
  (\email{zengzj22@mails.tsinghua.edu.cn}).}
 \and Junqing Chen\thanks{Department of Mathematical Sciences, Tsinghua University, Beijing, China  
  (\email{jqchen@tsinghua.edu.cn}).}
\and Zuoqiang Shi\thanks{Yau Mathematical Sciences Center, Tsinghua University, Beijing, China and Yanqi Lake Beijing Institute of Mathematical Sciences and Applications, Beijing, China  
  (\email{zqshi@tsinghua.edu.cn}).Corresponding author.}}

\usepackage{amsopn}

\nolinenumbers
\ifpdf
\hypersetup{
  pdftitle={Kernel Variational Inference Filter},
  pdfauthor={Weiye Gan, Zhijun Zeng, Junqing Chen, and Zuoqiang Shi}
}
\fi

\begin{document}

\maketitle

\begin{abstract}
    We present a novel particle flow for sampling called kernel variational inference flow (KVIF). KVIF do not require the explicit formula of the target distribution which is usually unknown in filtering problem. Therefore, it can be applied to construct filters with higher accuracy in the update stage. Such an improvement has theoretical assurance. Some numerical experiments for comparison with other classical filters are also demonstrated.
\end{abstract}

\begin{keywords}
nonlinear filtering problem, particle filter, variational inference
\end{keywords}

\section{Introduction}
In many fields, such as signal processing~\cite{roth2017ensemble}, robotics~\cite{chen2011kalman,thrun2002particle}, and economics~\cite{flury2011bayesian}, the task of estimating the state of a system over time, known as filtering, is fundamental. Filtering involves sequentially updating an estimate of the state of the system based on noisy observations and a dynamic model. Specifically, the model considered in this article is 

\begin{equation}\label{eq:nonlinear_filtering_problem}
\begin{array}{cccl}
     \text{State}:& x_{k+1} &=&f_k(x_{k}) + w_k  \\
     \text{Observation}:& y_k &=&h_k(x_k) + v_k 
\end{array},
\end{equation}
where $x_k\in \mathbb{R}^d$ is a vector of states in step $k$, $y_k\in \mathbb{R}^m$ is a vector of observations in step $k$ and $w_k$, $v_k$ are noises. The target of the filtering problem is to estimate the state $x_k$, given the values of the observations from the time step $0$ to $k$, that is, the values of $y_1,y_2,\dots,y_k$. 

When $f$ and $h$ are linear and the noises are Gaussian, the Kalman filter~\cite{welch1995introduction} can afford an analytic optimal estimate for $x_k$. However, in nonlinear filtering problems, obtaining an optimal estimate is usually infeasible. As a result, numerous methods have been developed in recent years to obtain suboptimal estimates for nonlinear systems, such as the Extended Kalman Filter (EKF)~\cite{ribeiro2004kalman}, Unscented Kalman Filter (UKF)~\cite{wan2000unscented}, ensemble Kalman filter (EnKF)~\cite{evensen2003ensemble} and Yau-Yau filter~\cite{yau2008real}. One promising technique that has emerged in the context of nonlinear filtering is particle filter (PF)~\cite{gordon1993novel}. Particle filter employs a set of random samples, or particles, to represent the posterior distribution $p(x_k|y_{1:k})$ and evolve over time according to the system dynamics and the measurement model. The core of the particle filter is the process of solving a Bayesian inference problem by updating particles based on new observations. The most basic version of particle filtering implements this through a resampling step, which is simple but comes with its own drawbacks. In particular, resampling can lead to issues such as sample degeneracy, where the majority of particles become ineffective, while large numbers of particles may be required to provide adequate approximations~\cite{carpenter1999improved}.

To address the challenges posed by particle filters, an alternative approach to solving the Bayesian inference problem is Variational Inference (VI)~\cite{blei2017variational}. Variational inference methods aim to approximate the true posterior distribution by optimizing a simpler, parameterized distribution. These methods are particularly useful when exact inference is intractable. Several VI techniques have been proposed, including diffusion posterior sampling methods~\cite{bao2024score,chu2025split} and particle-based variational inference~\cite{liu2016stein,chen2018unified,liu2019understanding,di2021neural,fan2021variational,zhang2024semi}. Diffusion posterior sampling methods utilize the diffusion model, which is a powerful tool for learning high-dimensional distributions. While the training may be computationally expensive. Particle-based variational inference, such as Stein Variational Gradient Descent (SVGD)~\cite{liu2016stein}, has attracted attention for its ability to transform a set of particles towards the posterior distribution by iteratively applying a gradient-based update.

SVGD relies on the score function, the gradient of the log-posterior, to update the particles. However, this requires an explicit analytic form of the score function, which makes it challenging to apply directly in particle filtering settings, where the posterior is typically unknown. In this paper, inspired by the strengths of both particle filtering and SVGD, we propose a novel filtering method called Kernel Variational Inference Flow Filter (KVIFF). This method combines the flexibility of particle filters with the efficiency of variational inference, offering a new approach to nonlinear filtering problem.

The paper is organized as follows. In \cref{sec:preliminary}, we introduce particle filter and SVGD which are closely related to our method. The complete process of KVIFF and some theoretical analysis are presented in \cref{sec:preliminary}. In \cref{sec:experiments}, we demonstrate some numerical results to compare our method with other classical filters. In \cref{sec:conclusions}, we conclude this paper and present several aspects for future research.

\section{Preliminary}
\label{sec:preliminary}
    \subsection{Particle Filter}
    Particle filter (PF)~\cite{gordon1993novel} is a sequential Monte Carlo\\
     method used to approximate the posterior distribution of the state \(x_k\) by representing it with a set of particles. Each particle \(x_k^{(i)}\) represents a sample of the posterior distribution $p(x_k|y_{1:k})$. The procedure proceeds in the following steps:
    
    1. Initialization: The initial state distribution \(p(x_0)\) is approximated by a set of particles \(\{x_0^{(i)}\}_{i=1}^N\) sampled from the known distrbution of the initial state.
       
       \begin{equation*}
           x_0^{(i)} \sim p(x_0), \quad i = 1, \dots, N
       \end{equation*}
    
    2. Prediction: Given the previous state estimate, the particles are propagated according to the state model \(f_k\), using a process noise model. The predicted state for each particle is:
    
       \begin{equation*}
           x_{k+1}^{(i)} = f_k(x_k^{(i)}) + w_k^{(i)}, \quad i = 1, \dots, N
       \end{equation*}
    where $w_k^{(i)}$ are samples of distribution of noise $w_k$.
    
    3. Update: Note that $x_k^{(i)}$ are the samples of $p(x_k|y_{1:k})$. $x_{k+1}^{(i)}$ can be considered as the samples of $p(x_{k+1}|y_{1:k})$. With Bayesian rule, we have
    \begin{equation}
        \label{eq:bayesian_rule}
        p(x_{k+1}|y_{1:{k+1}})=\frac{p(y_{k+1}|x_{k+1})p(x_{k+1}|y_{1:{k}})}{p(y_{k+1}|y_{1:{k}})}
    \end{equation}
    After obtaining the observation \(y_{k+1}\), the importance weights of each particle are given based on the likelihood of the observation given the predicted state as follows:
    
       \begin{equation*}
        \alpha_{k+1}^{(i)} = \frac{p(y_{k+1}|x_{k+1})}{p(y_{k+1}|y_{1:{k}})}
       \end{equation*}
    which can be calculated analytically if the explicit formula of the density function of $v_k$ is known. Resampling $x_{k+1}^{(i)}$ with weights $\alpha_{k+1}^{(i)}$ allows the particles to be updated using the observation information.
    
    4. Estimation: The estimate of the state $x_k$ is obtained by computing the mean of the particles:
    
       \[
       \hat{x}_k = \frac{1}{N}\sum_{i=1}^N x_k^{(i)}
       \]
    
    This process is iterated for each time step \(k = 1, 2, \dots\), leading to a sequence of state estimates \(\hat{x}_k\). 

    In particle filter, the update stage typically involves a resampling step to mitigate the issue of degeneracy, where most of the particles have negligible weights. However, this resampling process can introduce problems, such as the depletion of particle diversity, which may lead to a loss of information and performance degradation. To address these issues, several modifications to the standard particle filter have been proposed. One notable improvement is the Ensemble Kalman Filter (EnKF)~\cite{evensen2003ensemble}, which leverages a Kalman filter-like update to propagate the state estimates, while maintaining the diversity of the particle set. 

    \subsection{Stein Variational Gradient Descent}
    Variational inference try to find an optimal approximation of target distribution $p(x)$ in a predefined set $\mathcal{Q}$ by minimizing the KL divergence,
    \begin{equation*}
        q^* = \mathop{\arg\min}_{q\in\mathcal{Q}}KL(q|| p):=\mathop{\arg\min}_{q\in\mathcal{Q}}\mathbb{E}_q[\log q - \log p].
    \end{equation*}
    To strike a balance of accuracy and solvability, SVGD~\cite{liu2016stein} considers $\mathcal{Q}$ as sets consisting of distributions obtained by smooth transforms from a tractable initial distribution $q_0$. Particularly,
    \begin{equation*}
        \mathcal{Q}:=\{q_0^{[T]}|T(x)=x+\phi(x), \phi\in\mathcal{H}^d\}
    \end{equation*}
    where $q_0^{[T]}$ is the density of $z=T(x)$, $x\sim q_0$ and $\mathcal{H}^d$ is a RKHS with respect to a kernel $k$. In \cite{liu2016stein}, it has shown that
    \begin{equation*}
        \nabla_{\phi}KL(q_0^{[T]}|| p)|_{\phi=0}=-\mathbb{E}_{x\sim q_0}[k(x,\cdot)\nabla_x\log p(x)+\nabla_x k(x,\cdot)]
    \end{equation*}
    This suggests that $T(x):=x+\varepsilon\mathbb{E}_{x\sim q_0}[k(x,\cdot)\nabla_xp(x)+\nabla_x k(x,\cdot)]$ is equivalent to functional gradient descent in RKHS. Therefore, SVGD can approximate $p(x)$ with an iteration
    \begin{equation}
        \label{eq:SVGD}
        q_{t+1} = q_t + \varepsilon \mathbb{E}_{x\sim q_t}[k(x,\cdot)\nabla_x\log p(x)+\nabla_x k(x,\cdot)].
    \end{equation}
    This is also called continuous or mean-field SVGD. It has been shown that $q_t$ converges to the target distribution $p$ as $t\rightarrow\infty$ in the sense of KL divergence~\cite{liu2017stein}. This result is sharpen with weaker conditions or explicit rates in~\cite{korba2020non,chewi2020svgd,salim2022convergence,sun2023convergence,duncan2023geometry}. In practice, \cref{eq:SVGD} is usually approximated by a $n$-particles version,
    \begin{equation}
        \label{eq:SVGD_particles}
        x_{t+1}^{(i)}=x_{t}^{(i)}+\varepsilon\frac{1}{n}\sum_{j=1}^n\left\{k(x_{t}^{(j)},x_{t}^{(i)})\nabla\log p(x_{t}^{(j)})+\nabla k(x_{t}^{(j)},x_{t}^{(i)})\right\}
    \end{equation}
    for $i=1,\dots,n$. And the empiricle measure $\hat{q}_t(x):=\frac{1}{n}\sum_{i=1}^n\delta(x-x_t^{(i)})$ converges to $p$ at an order at least $1/\sqrt{\log\log n}$~\cite{shi2023finite}. This finite-particles convergence rate is improved double-exponentially to $1/\sqrt{n}$ in~\cite{banerjee2025improved}. 
    
    SVGD performs well when the number of samples is limited, and has been applied in other problems such as generative models~\cite{feng2017learning} and reinforcement learning~\cite{liu2017stein}. However, its performance highly depends on the choice of kernel $k$ and the quadratic computational complexity of the kernel matrix impedes the practical usage of a large number of particles. To construct more flexible gradient flow approximations, functional gradient flow methods~\cite{di2021neural,dong2022particle,cheng2023particle} can be adopted to expand the kernel space.
    
\section{Method}
\label{sec:method}
SVGD is a powerful tool to solve Bayesian inference problem. But it cannot be applied directly in \eqref{eq:bayesian_rule} since the score function of posterior distribution is unknown. However, we actually have some information about the posterior distribution. To simplify the notation, we rewrite \eqref{eq:bayesian_rule} as
\begin{equation*}
    p(x) = \frac{\tilde{Q}(x)q(x)}{C_Q}=Q(x)q(x)
\end{equation*}
where $p(x)=p(x_{k+1}|y_{1:k+1})$ is the target density function. Samples of $q(x)=p(x_{k+1}|y_{1:k})$ can be accessed. The explicit formula of $\tilde{Q}=p(y_{k+1}|x_{k+1})$ can be derived from the observation equation and $C_Q$ is the normalization constant which can be estimated with Monte-Carlo method. With such structure, one can slightly revise the update direction in \eqref{eq:SVGD} from
\begin{equation*}   
    \phi_t^{SVGD}:=\mathbb{E}_{x\sim q_t}[k(x,\cdot)\nabla_x\log p(x)]+\mathbb{E}_{x\sim q_t}[\nabla_x k(x,\cdot)].
\end{equation*}
to
\begin{equation*}   
    \phi_t:=\mathbb{E}_{x\sim p}[k(x,\cdot)\nabla_x\log p(x)]+\mathbb{E}_{x\sim q_t}[\nabla_x k(x,\cdot)].
\end{equation*}
With integration by parts,
\begin{equation*}
    \begin{aligned}
        &\quad\mathbb{E}_{x\sim p}[k(x,\cdot)\nabla_x\log p(x)]\\
        &=\mathbb{E}_{x\sim q}[k(x,\cdot)Q(x)(\nabla_x\log Q(x)+\nabla_x\log q(x))]\\
        &=\mathbb{E}_{x\sim q}[k(x,\cdot)\nabla_x Q(x)]-\mathbb{E}_{x\sim q}[\nabla_x(k(x,\cdot)Q(x))]\\
        &=-\mathbb{E}_{x\sim q}[\nabla_x k(x,\cdot) Q(x)].
    \end{aligned}
\end{equation*}
Hence,
\begin{equation}
    \label{eq:KVIF}
    \phi_t=-\mathbb{E}_{x\sim q}[\nabla_x k(x,\cdot) Q(x)]+\mathbb{E}_{x\sim q_t}[\nabla_x k(x,\cdot)]
\end{equation}
which can be calculated with the Monte Carlo method if we have samples of $q$ and $q_t$. This flow is called kernel variational inference flow (KVIF). Intuitively, if the difference between $q_t$ and $p$ is not significant, $\phi_t$ should be an approximation of $\phi_t^{SVGD}$ which is the steepest descent direction of the KL divergence. Therefore, $\phi_t$ should still be a descent direction. Rigorously, under a setting of continuous time, for SVGD, we have
\begin{equation*}
    \frac{\partial}{\partial t}q^{SVGD}(x,t)=-div(q^{SVGD}(x,t)\phi_t^{SVGD})
\end{equation*}
following the Fokker-Planck equation. Theorem 3.4 in~\cite{liu2017stein} shows that
\begin{equation*}
    \frac{d}{dt}KL(q_t^{SVGD}||p)=-\mathbb{D}(q_t^{SVGD}||p)^2
\end{equation*}
where $\mathbb{D}(q||p)$ is the Stein discrepancy. Since Stein discrepency is positive definite, SVGD can ensure a monotonically decreasing in the KL divergence and the convergence to the target distribution. For KVIF, we can demonstrate that it is a descent direction of weighted $L^2$ loss without any requirement on the distance between $q_t$ and $p$.
\begin{theorem}
\label{thm:mon-decrease}
    Let $k(x,x')$ be a symmetric kernel satisfying $\nabla_xk(x,x')=$\\
    $-\nabla_{x'}k(x,x')$ and $q(x,t)$ satisfy the Fokker-Planck equation
    \begin{equation*}
        \frac{\partial}{\partial t}q(x,t) = -div(q(x,t)\phi_t)=div\Big(q(x,t)\int_{\mathbb{R}^d}\nabla_s k(s,x)(p(s) - q(s,t))ds\Big)
    \end{equation*}
    where $\phi_t$ is defined as in \cref{eq:KVIF}. Then we have 
    \begin{equation*}
        \frac{d}{d t}\mathcal{L}(p,q(\cdot,t))=-2\int q(x,t)\Big(\int_{\mathbb{R}^d}\nabla_s k(s,x)(p(s) - q(s,t))ds\Big)^2dx\leq 0
    \end{equation*} 
    where $\mathcal{L}(p,q(\cdot,t))$ is the weighted $L^2$ loss defined as
    \begin{equation*}
        \mathcal{L}(p,q(\cdot,t)):=\int_{\mathbb{R}^d}\int_{\mathbb{R}^d}(p(x)-q(x))k(x,x')(p(x')-q(x'))dxdx'.
    \end{equation*}
\end{theorem}
\begin{proof}
    \begin{equation*}
        \begin{aligned}
            &\quad\frac{d}{d t}\mathcal{L}(p,q(\cdot,t))\\
            &=2\int_{\mathbb{R}^d}\int_{\mathbb{R}^d}\frac{\partial}{\partial t}q(x,t)k(x,x')(q(x',t)-p(x'))dxdx'\\
            &=-2\int q(x,t)\int_{\mathbb{R}^d}\nabla_s k(s,x)(p(s) - q(s,t))ds\int_{\mathbb{R}^d}\nabla_x k(x,x')(q(x',t)-p(x'))dx'dx\\
            &=-2\int q(x,t)\Big(\int_{\mathbb{R}^d}\nabla_s k(s,x)(p(s) - q(s,t))ds\Big)^2dx\\
            &\leq 0
        \end{aligned}
    \end{equation*}
    where we applied integration by parts and assumptions of kernel $k$ in the last two equations respectively.
\end{proof}
If kernel $k(x,x')$ is positive definite, with the Mercer's theorem, we have the spectral decomposition
\begin{equation*}
    k(x,x')=\sum_{j=1}^\infty \lambda_je_j(x)e_j(x')
\end{equation*}
where $\lambda_1\geq\lambda_2\geq\dots>0$ are the eigenvalues and $\{e_j\}_{j=1}^\infty$ are the eigenfunctions which form an orthogonal basis of the $L^2(\mathbb{R}^d)$ space. We can define a Hilbert space 
\begin{equation*}
    \mathcal{H}=\Bigg\{f=\sum_{j=1}^\infty f_je_j\Bigg|\sum_{j=1}^\infty \lambda_jf_j^2<\infty\Bigg\}
\end{equation*}
equipped with inner product
\begin{equation*}
    \langle f,g\rangle_{\mathcal{H}}=\int_{\mathbb{R^d}}\int_{\mathbb{R^d}}f(x)k(x,x')g(x')dxdx'=\sum_{j=1}^\infty \lambda_jf_jg_j.
\end{equation*}
Then $\mathcal{L}(p,q(\cdot,t))=\lVert p-q(\cdot,t)\rVert_\mathcal{H}^2\geq 0$. And \cref{thm:mon-decrease} demonstrates that $\lVert p-q(\cdot,t)\rVert_\mathcal{H}^2$ monotonically decreases as $t$ increases. Hence, we have $\frac{d}{dt}\mathcal{L}(p,q(\cdot,t))\rightarrow 0$ as $t\rightarrow\infty$. And the following lemma shows that under some regular assumptions on $p$ and $q$, $\frac{d}{dt}\mathcal{L}(p,q)=0$ if and only if $p=q$ almost surely.
\begin{lemma}
    Suppose that kernel $k$ is positive definite and satisfies the assumptions in \cref{thm:mon-decrease}. $p,q$ are density functions of class $C^1\cap L^2$ and satisfies
    \begin{equation*}
        \int_{\mathbb{R}^d}\nabla_x(k(x,x')(p(x)-q(x)))dx=0.
    \end{equation*}
    Then, $\frac{d}{dt}\mathcal{L}(p,q)=0$ if and only if $p=q$ almost surely.
\end{lemma}
\begin{proof}
    Denote the support of $q$ by $S=\bigcup S_i$ where $\{S_i\}$ are connected components.  Note that $q\geq 0$. With $\frac{d}{dt}\mathcal{L}(p,q)=0$, we have
    \begin{equation*}
        0=\int_{\mathbb{R}^d}\nabla_s k(s,x)(p(s) - q(s))ds=-\int_{\mathbb{R}^d} k(s,x)(\nabla p(s) - \nabla q(s))ds
    \end{equation*}
    almost surely on $S$. Since $k$ is positive definite, we have
    \begin{equation*}
        \nabla p(x) = \nabla q(x)
    \end{equation*}
    almost surely on $S$, which means that there exists constants $\{C_i\}$ such that 
    \begin{equation*}
        q(x) = p(x) + C_i
    \end{equation*}
    for $x\in S_i$ almost surely. Since $q$ is continuous, we have $p(x)=0$ for $x\in \partial S_i$. Therefore, to ensure that $q\geq 0$, $C_i$ is nonnegative. Hence,
    \begin{equation*}
        1=\sum_{i}\int_{S_i}p(x)dx=\sum_{i}\int_{S_i}q(x)-C_idx\leq \sum_{i}\int_{S_i}q(x)dx\leq 1
    \end{equation*}
    which derives that $C_i=0$ and $q(x)=0$ almost surely on $S^c$. This means that $p=q$ almost surely.
\end{proof}

For model \eqref{eq:nonlinear_filtering_problem}, one can apply particle filter by solving the Bayesian posterior problem \eqref{eq:bayesian_rule} with KVIF. This filter is called kernel variational inference flow filter (KVIFF), which can be concluded as in \cref{alg:kvif}.
\begin{remark}
\label{remark1}
For the selection of $\{x_{k,0}^{(i)}\}_{i=1}^N$ in row 10 of \cref{alg:kvif}, a direct way is to set $x_{k,0}^{(i)}=x_k^{(i)}$. However, to improve efficiency, a more recommended approach is to let $\{x_{k,0}^{(i)}\}_{i=1}^N$ be the updated particles with other simpler filter methods, such as particle filter and ensemble Kalman filter. With such a choice, KVIFF can also be understood as PF or EnKF with additional iterations for variational inference to enhance accuracy in the update stage.
\end{remark}
\begin{algorithm}[h]
\caption{Kernel Variational Inference Flow Filter}
\label{alg:kvif}
\begin{algorithmic}[1]
\REQUIRE Initial particals $\{x_0^{(i)}\}_{i=1}^N$, observation $\{y_k\}_{k=1}^K$, kernel $k$, step size $\varepsilon$, number of steps $N_s$
\ENSURE particals $\{x_k^{(i)}\}_{i=1}^N$ and estimates $\hat{x}_k$ for $k=1,\dots , K$

\FOR{$k=1$ to $K$}
    \STATE \textbf{Prediction:}
    \FOR{$i=1$ to N}
        \STATE sample $w_k^{(i)}\sim w_k$
        \STATE $x_{k}^{(i)}\gets f_k(x_{k-1}^{(i)}) + w_k^{(i)}$
    \ENDFOR
    \STATE \textbf{Update:}
    \STATE Compute likelihood function $\tilde{Q}$ of $p(y_k|x_k)$ via $y_k$ and the measurement equation in \eqref{eq:nonlinear_filtering_problem}
    \STATE $C_Q=\frac{1}{N}\sum_{i=1}^N \tilde{Q}(x_{k}^{(i)})$
    \STATE Select $\{x_{k,0}^{(i)}\}_{i=1}^N$
    \FOR{$\tau=1$ to $N_S$}
        \FOR{$i=1$ to N}
            \STATE $x_{k,\tau}^{(i)}=x_{k,\tau-1}^{(i)}+\frac{\varepsilon}{N}\sum_{j=1}^N\{\nabla k(x_{k}^{(j)},x_{k,\tau}^{(i)})\tilde{Q}(x_{k}^{(j)})/C_Q-\nabla k(x_{k,\tau-1}^{(j)},x_{k,\tau}^{(i)})\}$
        \ENDFOR
    \ENDFOR
    \STATE $x_k^{(i)}\gets x_{k,N_s}^{(i)}$
    \STATE $\hat{x}_k\gets \frac{1}{N}\sum_{i=1}^N x_k^{(i)}$
\ENDFOR
\end{algorithmic}
\end{algorithm}

\section{Experimental results}
\label{sec:experiments}
In this section, we present several numerical experiments to compare KVIFF with other classical filters including particle filter (PF), Kalman filter (KF), ensemble Kalman filter (EnKF) and Yau-Yau filter by computing the $L^2$ error between the estimated trajectories and the ground truth. All experiments are repeated 10 times, and the median results are shown to avoid the influence of random seeds. Kalman filter is only tested in linear filter problem with Gaussian noises. Yau-Yau filter is implemented following the method in \cite{yueh2014efficient} and is only tested in low-dimensional cases. As mentioned in \cref{remark1}, at each time step, we update the particles with PF or EnKF before the iteration of KVIFF. Hence, it is significant to observe whether KVIFF performs better than these two methods. All computations are done using a NVIDIA GeForce RTX 3060 Laptop GPU.
\subsection{Linear filter problem} 
\begin{figure}[htbp]
    \centering
    \begin{subfigure}{1.0\textwidth}
        \includegraphics[width=\linewidth]{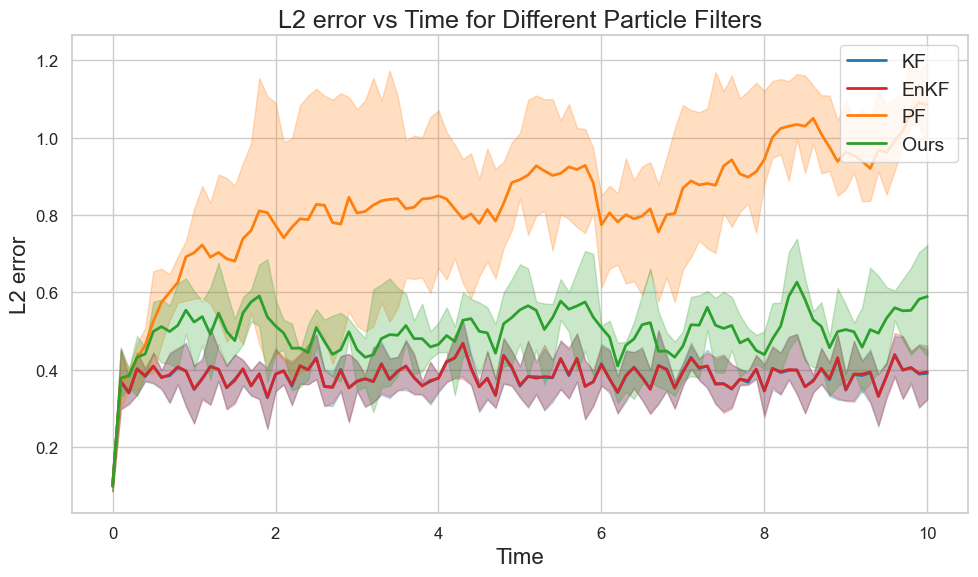}
        \caption{$L^2$ error}
    \end{subfigure}
    
    \vspace{\floatsep} 
    
    \begin{subfigure}{1.0\textwidth}
        \includegraphics[width=\linewidth]{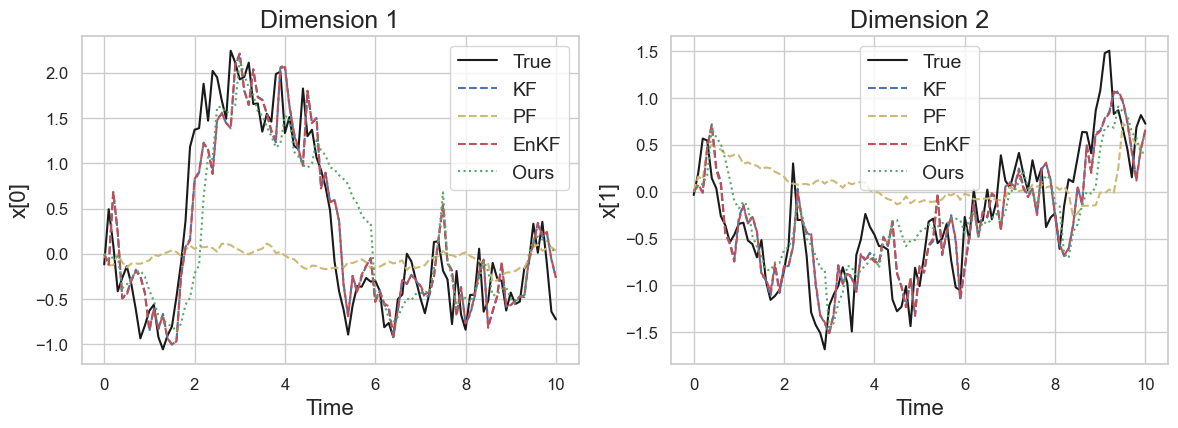}
        \caption{trajectories}
    \end{subfigure}
    \caption{$L^2$ error and trajectories (first two dimensions) for linear filter problem}
    \label{figure:linear}
\end{figure}

\begin{figure}[htbp]
    \centering
    \begin{subfigure}{1.0\textwidth}
        \includegraphics[width=\linewidth]{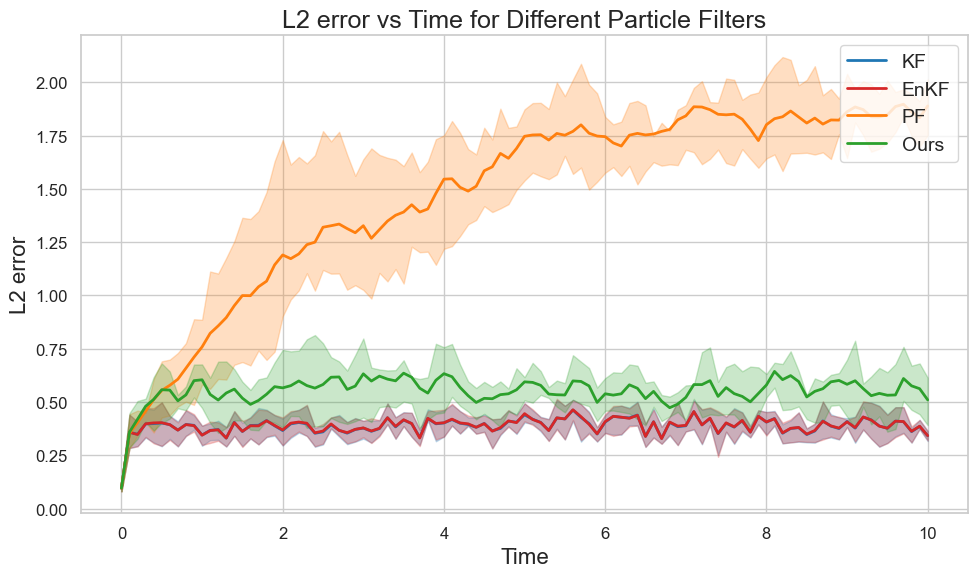}
        \caption{$L^2$ error}
    \end{subfigure}
    
    \vspace{\floatsep} 
    
    \begin{subfigure}{1.0\textwidth}
        \includegraphics[width=\linewidth]{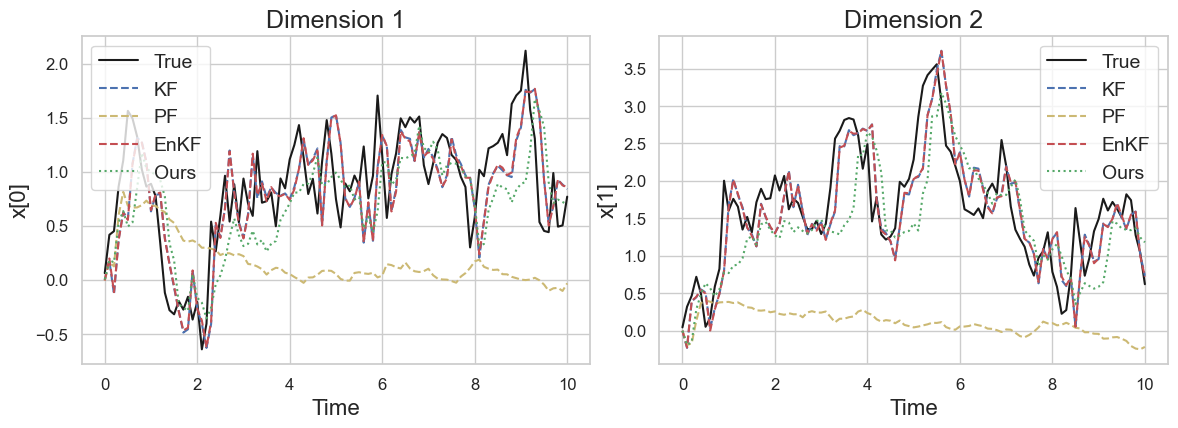}
        \caption{trajectories}
    \end{subfigure}
    \caption{$L^2$ error and trajectories (first two dimensions) for linear filter problem with bias}
    \label{figure:linear_bias}
\end{figure}

\begin{figure}[htbp]
    \centering
    \begin{subfigure}{1.0\textwidth}
        \includegraphics[width=\linewidth]{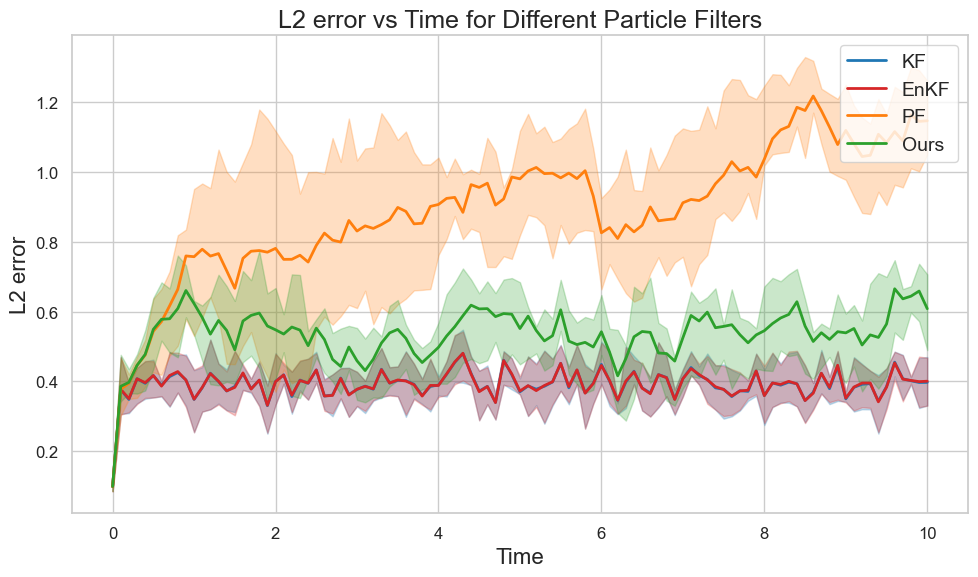}
        \caption{$L^2$ error}
    \end{subfigure}
    
    \vspace{\floatsep} 
    
    \begin{subfigure}{1.0\textwidth}
        \includegraphics[width=\linewidth]{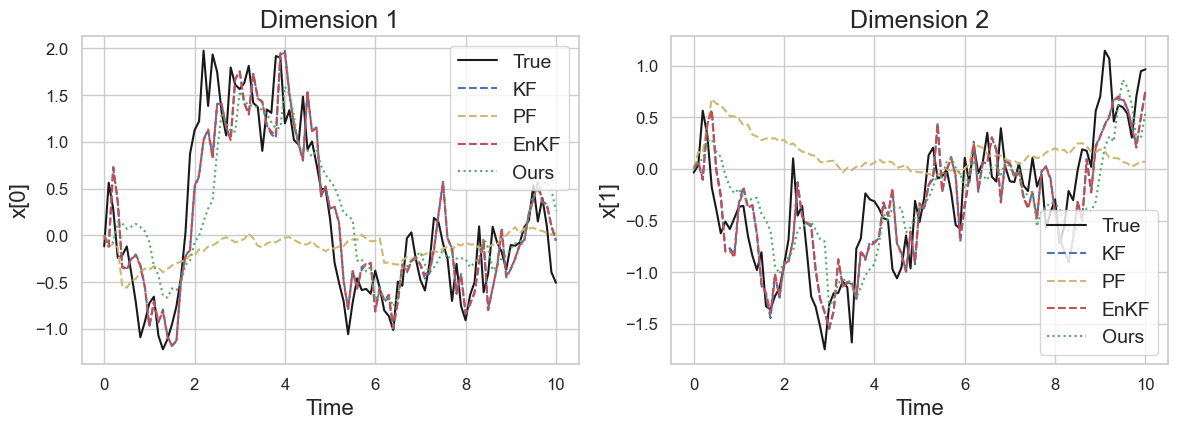}
        \caption{trajectories}
    \end{subfigure}
    \caption{$L^2$ error and trajectories (first two dimensions) for linear filter problem with correlated noise}
    \label{figure:linear_cor}
\end{figure}

Our first numerical experiment is a 10D linear filter problem.
\begin{equation}
    \label{eq:linear_filter}
\begin{cases}
x_{k+1} = f(x_k) + w_k,\\[4pt]
y_k = x_k + v_k,
\end{cases}
\end{equation}
with
\[
f(x) = x+\Delta tAx,
\quad
A = [a_{ij}],\;
a_{ij} = \begin{cases}
0.1, & i+1 = j,\\
-0.5,& i=j,\\
0,   & \text{otherwise},
\end{cases}
\]
and $\Delta t=0.1$, $\mathbb{E}[w_kw_k^T]=\Delta tI$, $\mathbb{E}[v_kv_k^T]=0.1I$. Time step $K=100$ so the total time $T=K\Delta t=10$. All methods are implemented with 1000 particles. For the setting of KVIFF, $\{x_{k,0}^{(i)}\}_{i=1}^N$ are selected as the updated particles with PF. Note that with this selection, KVIFF can be considered as PF with additional iteration in the update stage. Kernel $k(x,x')=\exp(\lVert x-x'\rVert^2/10)$. Step size $\varepsilon=1e-3$ and the number of steps $N_s=50$. The results are shown in \cref{figure:linear}. Since this is a linear filter problem with Gaussian noise, KF offers the optimal estimation. And EnKF approximates KF well in this linear case. Our method has similar accuracy to KF and performs much better than PF. This verifies the effectiveness of the KVIF iteration. Moreover, when we consider the case where the modeling is not accurate, our method still shows a strong calibration ability. Specifically, when generating states and observations, we set the noise $w_k$ with non-zero mean $\mathbb{E}[w_k]=(0.2,\dots,0.2)^T$. We also suppose that we do not know such a change when constructing filters. Hence, all filter methods are still implemented with model \eqref{eq:linear_filter}. The results are shown in \cref{figure:linear_bias}. One can see that PF is severely affected by modeling bias. But our method can adapt this bias and maintain a competitive effect with KF and EnKF. We also test the case with actual correlated noise, i.e. $\mathbb{E}[w_kw_k^T]=\Delta t(I+E)$ where 
\[E_{ij}=\begin{cases}
0.3, & i+1 = j,\\
0,   & \text{otherwise}.
\end{cases}\]
The results are shown in \cref{figure:linear_cor}. Our method obviously outperforms PF.

\subsection{Cubic sensor problem} 
\begin{figure}[htbp]
    \centering
    \begin{subfigure}{1.0\textwidth}
        \includegraphics[width=\linewidth]{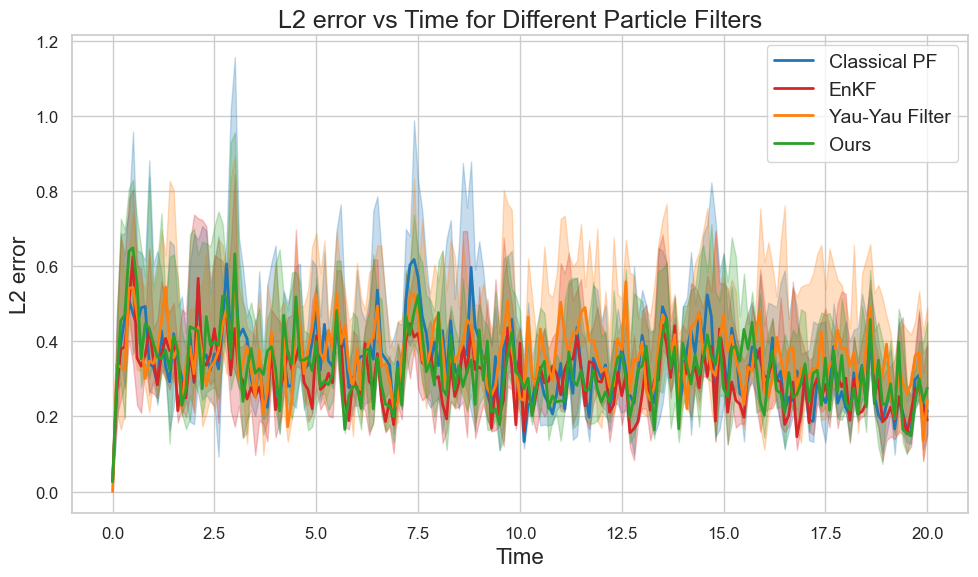}
        \caption{$L^2$ error}
    \end{subfigure}
    
    \vspace{\floatsep} 
    
    \begin{subfigure}{1.0\textwidth}
        \includegraphics[width=\linewidth]{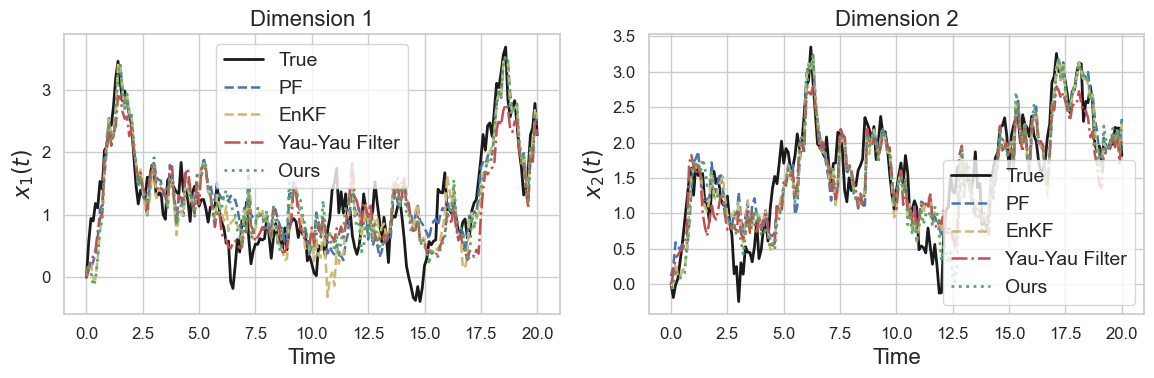}
        \caption{trajectories}
    \end{subfigure}
    \caption{$L^2$ error and trajectories for low dimensional cubic sensor problem}
    \label{figure:cubic_low_dimension}
\end{figure}

\begin{figure}[htbp]
    \centering
    \begin{subfigure}{1.0\textwidth}
        \includegraphics[width=\linewidth]{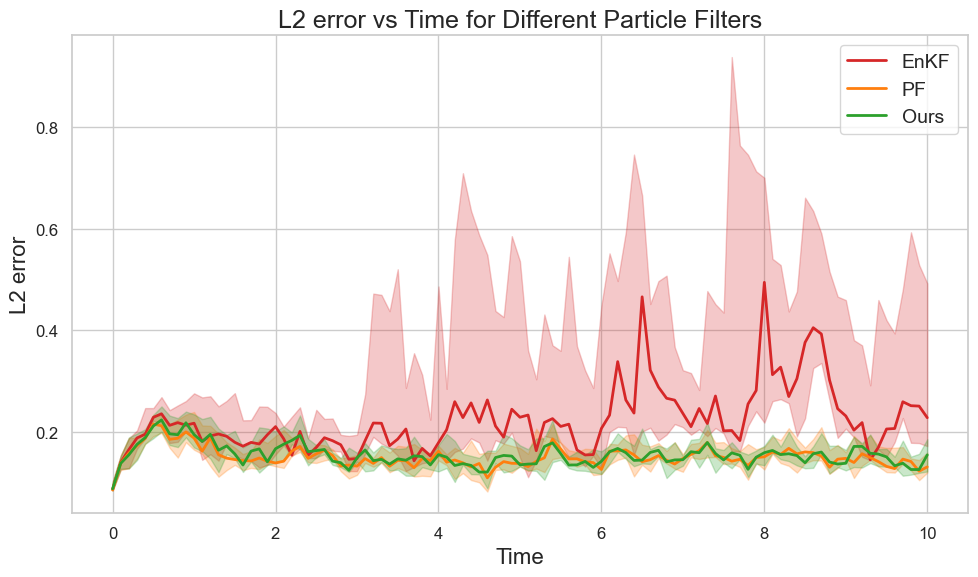}
        \caption{$L^2$ error}
    \end{subfigure}
    
    \vspace{\floatsep} 
    
    \begin{subfigure}{1.0\textwidth}
        \includegraphics[width=\linewidth]{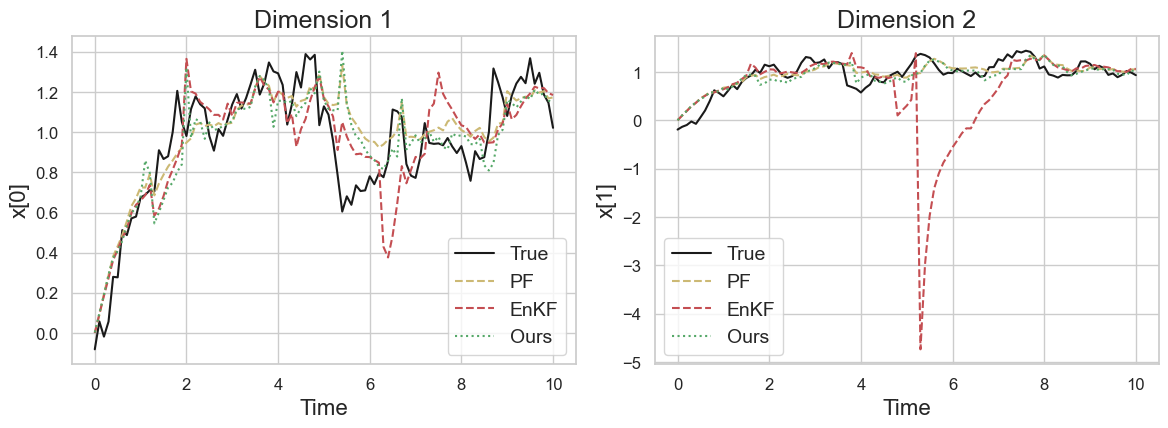}
        \caption{trajectories}
    \end{subfigure}
    \caption{$L^2$ error and trajectories (first two dimensions) for cubic sensor problem with heavy-tailed noise}
    \label{figure:cubic_heavy_tailed_noise}
\end{figure}

\begin{figure}[htbp]
    \centering
    \begin{subfigure}{1.0\textwidth}
        \includegraphics[width=\linewidth]{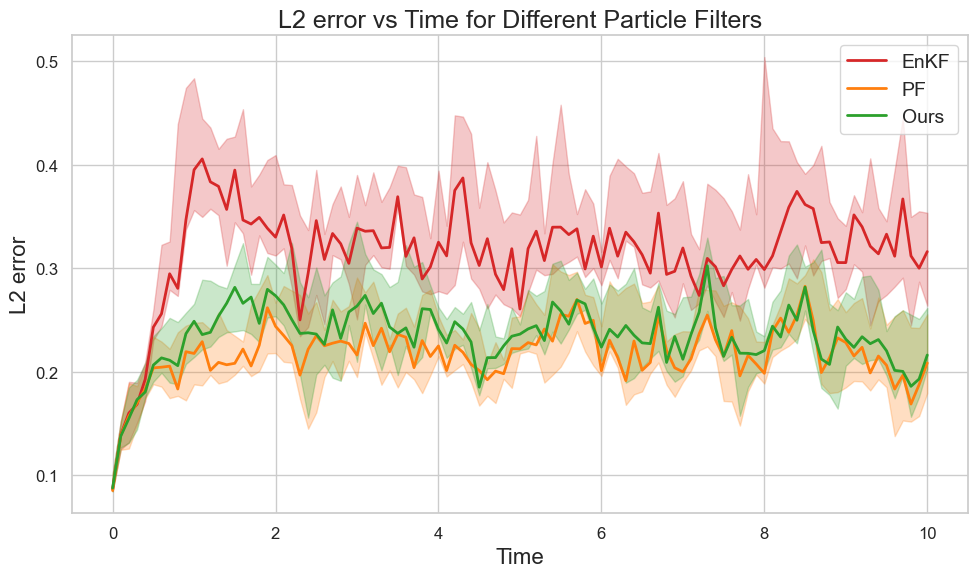}
        \caption{$L^2$ error}
    \end{subfigure}
    
    \vspace{\floatsep} 
    
    \begin{subfigure}{1.0\textwidth}
        \includegraphics[width=\linewidth]{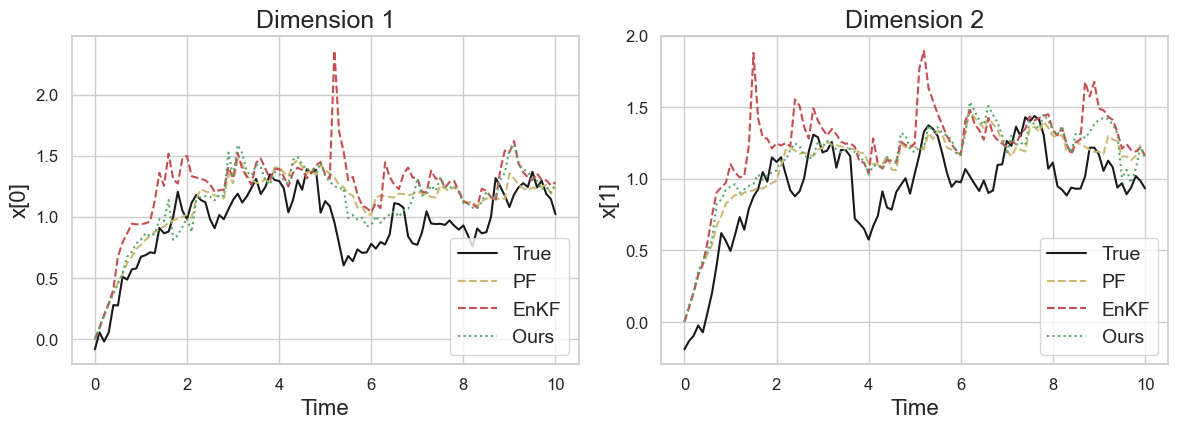}
        \caption{trajectories}
    \end{subfigure}
    \caption{$L^2$ error and trajectories (first two dimensions) for cubic sensor problem with asymmetric noise and bias}
    \label{figure:cubic_asymmetric_noise}
\end{figure}

The next problem is the cubic sensor problem which is focused by a series of research on Yau-Yau filter. To compare with Yau-Yau filter, we first consider a low-dimensional (2D) case studied in \cite{yueh2014efficient}.
\[
\begin{cases}
x_{k+1} = x_k + \Delta t\cos(x_k) + w_k,\\[4pt]
y_k = y_{k-1} + \Delta t x_k^3 + v_k.
\end{cases}
\]
and $\Delta t=0.1$, $\mathbb{E}[w_kw_k^T]=\Delta tI$, $\mathbb{E}[v_kv_k^T]=\Delta tI$. Time step $K=200$ so the total time $T=K\Delta t=20$. All methods except the Yau-Yau filter are implemented with 50 particles. For the setting of KVIFF, $\{x_{k,0}^{(i)}\}_{i=1}^N$ are still selected as the updated particles with PF. Kernel $k(x,x')=\exp(\lVert x-x'\rVert^2)$. Step size $\varepsilon=1e-3$ and the number of steps $N_s=50$. The results are shown in \cref{figure:cubic_low_dimension}, our method performs comparably with several other methods. However, for high-dimensional cases, the Yau-Yau filter faces an extremely high computational cost since it needs to solve high-dimensional PDEs. Moreover, it cannot handle non-Gaussian noises from the theoretical aspect. This is also a weakness for EnKF. While our method can flexibly deal with various noise distributions. To verify this point, we generate 10D states and observations with the measurement noise $v_k$ following a Cauchy distribution and implement the filter methods with the assumption that $v_k$ follows normal distribution. All methods are implemented with 1000 particles. Kernel $k(x,x')=\exp(\lVert x-x'\rVert^2/10)$. Other setting of KVIFF is the same as the low-dimensional case. The results are shown in \cref{figure:cubic_heavy_tailed_noise}. In this case, ensemble Kalman filter fails while our method and particle filter still work. We also apply log-normal distribution which is an asymmetric distribution to replace the normal distribution for $v_k$. And we suppose that there is a modeling bias for $w_k$, i.e. $\mathbb{E}[w_k]=(0.3,\dots,0.3)^T$ when generating data and $\mathbb{E}[w_k]=\mathbf{0}$ in the model for filter. The results are shown in \cref{figure:cubic_asymmetric_noise}. Our method is significantly better than particle filter and ensemble Kalman filter.

\subsection{Multi-target tracking problem} 
\begin{figure}[htbp]
    \centering
    \begin{subfigure}{1.0\textwidth}
        \includegraphics[width=\linewidth]{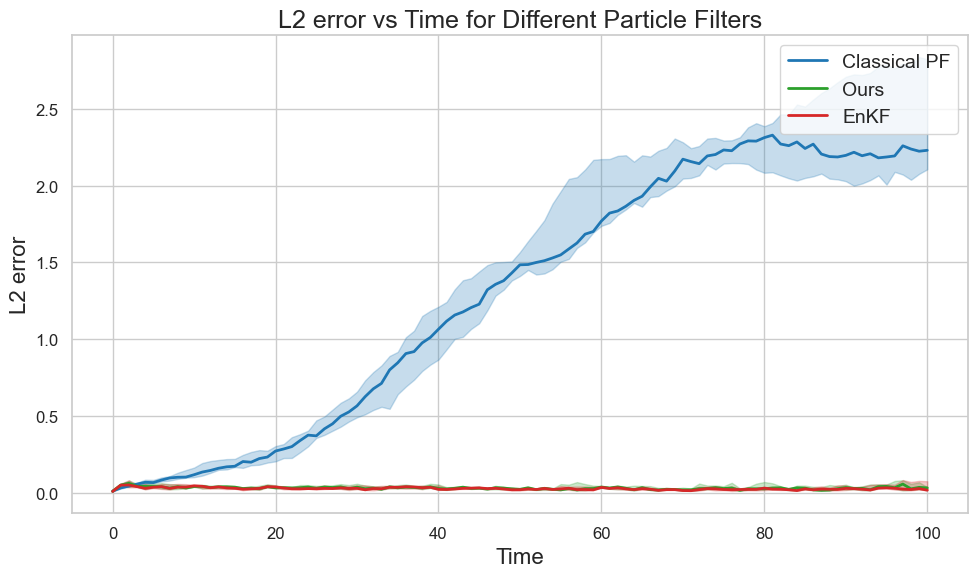}
        \caption{$L^2$ error}
    \end{subfigure}
    
    \vspace{\floatsep} 
    
    \begin{subfigure}{1.0\textwidth}
        \includegraphics[width=\linewidth]{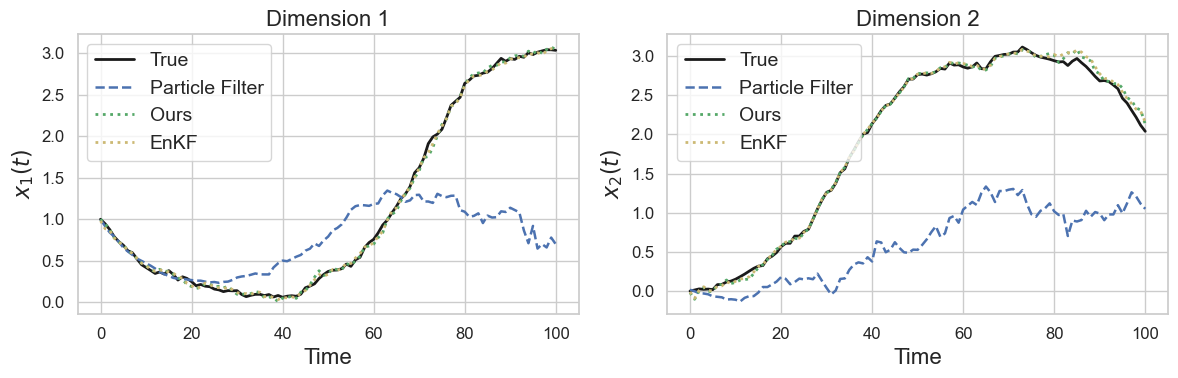}
        \caption{trajectories}
    \end{subfigure}
    \caption{$L^2$ error and the trajectory of the first object for multi-target tracking problem}
    \label{figure:multi-target_tracking}
\end{figure}

\begin{figure}[htbp]
    \centering
    \begin{subfigure}{1.0\textwidth}
        \includegraphics[width=\linewidth]{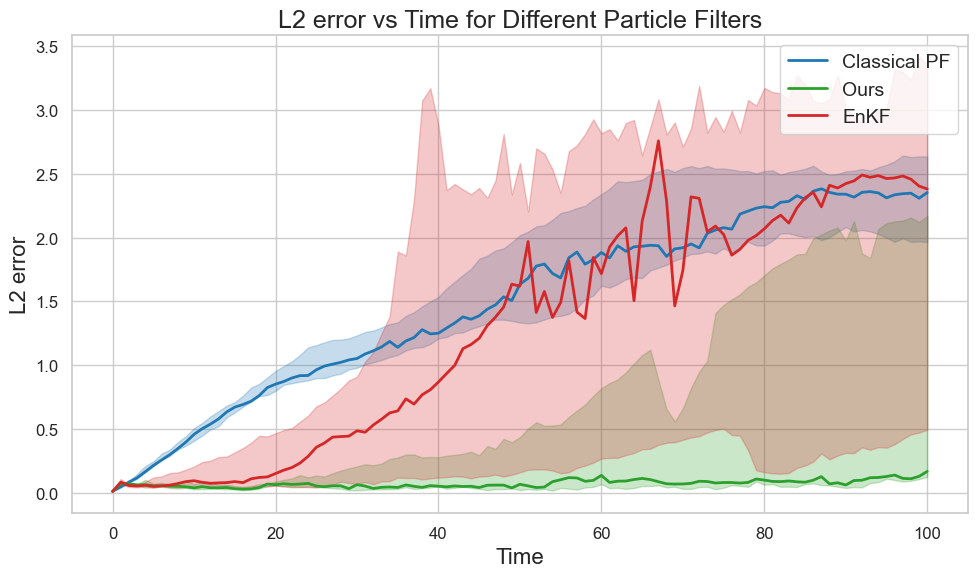}
        \caption{$L^2$ error}
    \end{subfigure}

    \vspace{\floatsep} 
    
    \begin{subfigure}{1.0\textwidth}
        \includegraphics[width=\linewidth]{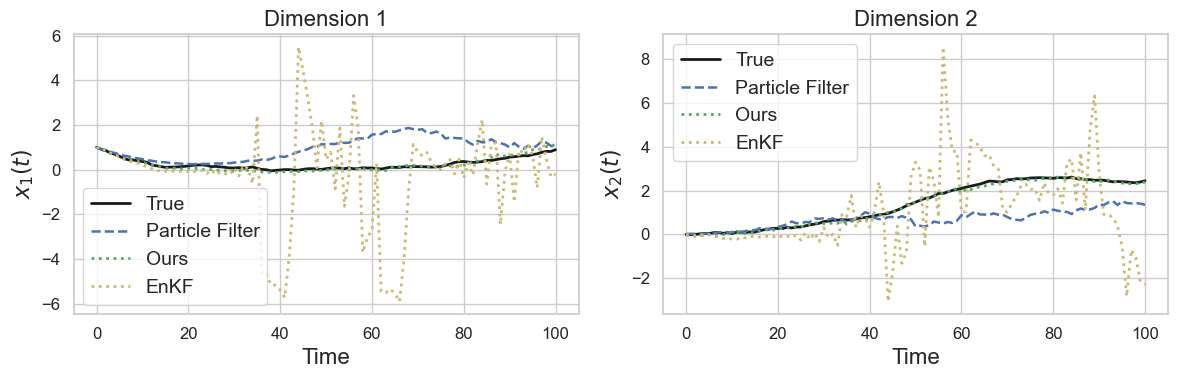}
        \caption{trajectories}
    \end{subfigure}
    \caption{$L^2$ error and the trajectory of the first object for multi-target tracking problem with bias}
    \label{figure:multi-target_tracking_bias}
\end{figure}

Finally, we consider multi-target tracking problem which is a more practical problem also focused in \cite{hlinka2011distributed,ding2012implementation}. For the dynamic system, we consider $M=4$ objects moving independently in a 2D velocity field with initial position $(1,0)$, $(2,0)$, $(3,0)$ and $(4,0)$.
\begin{equation}
    \label{eq:multi_target_state}
    x_{k+1}^{(p)} = x_k^{(p)} + v(x_k^{(p)}) \Delta t + w_k^{(p)}
\end{equation}
where $p=1,\dots,M$ where $\Delta t=0.1$, $\mathbb{E}[w_kw_k^T]=0.1\Delta tI$, $v$ is the Hamiltonian vector field associated with $H(x)=\sin(x_1)\sin(x_2)$, i.e. 
\[
v(x)=\Big(-\frac{\partial H}{\partial x_2}, \frac{\partial H}{\partial x_1}\Big)=(-\sin(x_1)\cos(x_2),\cos(x_1)\sin(x_2))
\] 
For the measurement equation, we model a wireless sensor network consisting of 25 acoustic amplitude sensor nodes deployed on a uniform grid of region $[-4,4]\times[-4,4]$. The measurement $y_k^j$ obtain by sensor $j$ at time $k$ is 
\begin{equation*}
    y_k^j=\sum_{p=1}^M\frac{A}{\lVert x_k^{(p)}-\xi_j\rVert_2+d_0} + v_k^j.
\end{equation*}
$\mathbb{E}[v_kv_k^T]=0.01\Delta tI$. $A=10$ represents the amplitude emitted by the objects. $d_0=0.1$ is a threshold that determines the maximum measurable amplitude. $\xi_j$ is the location of the sensor $j$. Time step $K=100$ so the total time $T=K\Delta t=10$. All methods are implemented with 500 particles. For the setting of KVIFF, $\{x_{k,0}^{(i)}\}_{i=1}^N$ are selected as the updated particles with EnKF since PF performs poorly on this task, see in \cref{figure:multi-target_tracking}. Kernel $k(x,x')=\exp(\lVert x-x'\rVert^2/10)$. Step size $\varepsilon=5e-5$ and the number of steps $N_s=200$. The results are shown in \cref{figure:multi-target_tracking}. The performance of EnKF is clearly better than that of particle filter. Our method maintains a comparable accuracy. For the next, we additionally consider a modeling bias for the velocity field. Specifically, the real states for generating data satisfy the following equation,
\begin{equation*}
    x_{k+1}^{(p)} = x_k^{(p)} + v(x_k^{(p)}) \Delta t - x_k^{(p)}\Delta t/5 + w_k^{(p)}
\end{equation*}
but we still adopt \eqref{eq:multi_target_state} in the filtering methods. We neglect a component of the velocity that points to the origin. The results are shown in \cref{figure:multi-target_tracking_bias}. KVIFF demonstrates significantly stronger calibration capability compared to other methods.

\section{Conclusions}
\label{sec:conclusions}
In this paper, we present a novel particle based filter method called Kernel Variational Inference Flow Filter (KVIFF). In the update stage, a VI-liked iteration is introduced to enhance accuracy with theoretical guarantee. Our method can be easily combined with other filters by selecting different initial particles for the iteration. The numerical experiments verify the improvement for KVIFF comparing to PF and EnKF, especially when the modeling of filter problem is not accurate. There are also some limitations and several potential directions to extend. Note that KVIFF is constructed based on a kernel liked SVGD. The issues faced by SVGD also need to be resolved for KVIFF. The performance of our method may be highly dependent on the choice of the kernels. The computational cost of KVIFF mainly comes from the quadratic computational complexity of the kernel matrix which also hinders the pratical usage of a large number of particles. To overcome these shortcomings, some following work of SVGD called functional gradient flow methods~\cite{di2021neural,dong2022particle,cheng2023particle} may be helpful. We leave this to the future work.


\bibliographystyle{siamplain}
\bibliography{references}
\end{document}